\documentclass[leqno,draft]{article}

\begin{document}

\title{Some remarks about Clifford analysis \\ and fractal sets}

\author{Stephen Semmes \\
        Rice University}

\date{}

\maketitle

        Let $n$ be a positive integer, and let $\mathcal{C}(n)$ be the
real \emph{Clifford algebra} with $n$ generators $e_1, \ldots, e_n$.
Thus $\mathcal{C}(n)$ is an associative algebra over the real numbers
with multiplicative identity element $1$ such that
\begin{equation}
        e_1^2 = e_2^2 = \cdots = e_n^2 = -1
\end{equation}
and
\begin{equation}
        e_j \, e_l = - e_l \, e_j \quad\hbox{when } j \ne l.
\end{equation}
This reduces to the complex numbers ${\bf C}$ when $n = 1$, and to the
quaternions ${\bf H}$ when $n = 2$.  For each $n$, $\mathcal{C}(n)$
has dimension $2^n$ as a vector space over the real numbers.

        Consider the first-order differential operators
$\mathcal{D}_L$, $\mathcal{D}_R$ acting on $\mathcal{C}(n)$-valued
functions on ${\bf R}^n$ defined by
\begin{equation}
  \mathcal{D}_L f = \sum_{j = 1}^n e_j \, \frac{\partial f}{\partial x_j}
        \quad\hbox{and}\quad
  \mathcal{D}_R f = \sum_{j = 1}^n \frac{\partial f}{\partial x_j} \, e_j.
\end{equation}
Of course, $\mathcal{D}_L f = \mathcal{D}_R f$ when $f$ is
real-valued, since real numbers commute with $e_1, \ldots, e_n$ by
definition.  Also,
\begin{equation}
        \mathcal{D}_L^2 f = \mathcal{D}_R^2 f
           = - \sum_{j = 1}^n \frac{\partial^2 f}{\partial x_j^2}.
\end{equation}
We say that $f$ is left or right \emph{Clifford holomorphic} on an
open set $U \subseteq {\bf R}^n$ if $\mathcal{D}_L f = 0$ or
$\mathcal{D}_R f = 0$ on $U$, respectively.

        There are variants of these notions for functions on ${\bf
R}^{n + 1}$, in which the derivative in the extra dimension does not
have a coefficient.  However, in analogy with the $\partial$ and
$\overline{\partial}$ operators from complex analysis, one considers
pairs of operators on the left and right, with and without an
additional minus sign for the derivatives in the other $n$ dimensions.
The product of two such operators again reduces to the Laplacian.
There are also variants for the quaternions with a derivative in
another direction, corresponding to the third imaginary basis vector
in ${\bf H}$.

        For the sake of simplicity, let us restrict our attention to
the version of Clifford analysis on ${\bf R}^n$ determined by the
operators $\mathcal{D}_L$, $\mathcal{D}_R$.  If $h$ is a harmonic
function on an open set in ${\bf R}^n$, then $\mathcal{D}_L h$,
$\mathcal{D}_R h$ are left and right Clifford holomorphic on the same
open set, respectively.  These functions are the same when $h$ is
real-valued, and hence both left and right Clifford holomorphic.  This
is basically the gradient of $h$, expressed as a Clifford-valued
function.  As in one complex variable, the components of Clifford
holomorphic functions are harmonic.

        The sum of two left or two right Clifford holomorphic
functions is left or right Clifford holomorphic, respectively.
Because of the noncommutativity of the Clifford algebra when $n \ge
2$, the product of two Clifford-holomorphic functions is not
necessarily Clifford holomorphic.  If $a \in \mathcal{C}(n)$ and $f$
is left Clifford-holomorphic, then $f \, a$ is left Clifford
holomorphic, and similarly $a \, f$ is right Clifford holomorphic when
$f$ is.  The product of a constant and a Clifford holomorphic function
in the other order is not necessarily Clifford holomorphic.

        These properties of products are also reflected in the
identities
\begin{equation}
 \mathcal{D}_L (f_1 \, f_2) = (\mathcal{D}_L f_1) \, f_2
  + \sum_{j = 1}^n e_j \, f_1 \, \frac{\partial f_2}{\partial x_j}
\end{equation}
and
\begin{equation}
 \mathcal{D}_R (f_1 \, f_2) =
  \sum_{j = 1}^n \frac{\partial f_1}{\partial x_j} \, f_2 \, e_j
     + f_1 \, (\mathcal{D}_R f_2).
\end{equation}

        The function
\begin{equation}
        E(x) = \frac{\sum_{j = 1}^n x_j \, e_j}{|x|^n}
\end{equation}
is left and right Clifford holomorphic for $x \ne 0$ and each $n$, and
in fact is a nonzero real multiple of the fundamental solution for
$\mathcal{D}_L$ and $\mathcal{D}_R$.  This means that $\mathcal{D}_L E
= \mathcal{D}_R E$ is a nonzero real multiple of the delta function at
$0$ in the sense of distributions.  This is because $E$ is a nonzero
real multiple of $\mathcal{D}_L h = \mathcal{D}_R h$, where $h(x) =
|x|$ when $n = 1$, $h(x) = \log |x|$ when $n = 2$, and $h(x) = |x|^{2
- n}$ when $n \ge 3$, and $h$ is a multiple of the fundamental
solution for the Laplacian in these dimensions.  If $f$ is a
continuously-differentiable $\mathcal{C}(n)$-valued function on ${\bf
R}^n$ with compact support, for instance, then it follows that
\begin{equation}
  f(x) = c_n \, \int_{{\bf R}^n} E(x - y) \, \mathcal{D}_L f(y) \, dy
       = c_n \, \int_{{\bf R}^n} \mathcal{D}_R f(y) \, E(x - y) \, dy
\end{equation}
for some real number $c_n \ne 0$ depending only on $n$ and every $x
\in {\bf R}^n$.

        As a nice little consequence of this formula, every continuous
$\mathcal{C}(n)$-valued function on a compact set $A \subseteq {\bf
R}^n$ with Lebesgue measure $0$ can be approximated uniformly by the
restrictions of left or right Clifford holomorphic functions on
neighborhoods of $A$.  For such a continuous function can be
approximated uniformly first by the restrictions to $A$ of
continuously-differentiable functions on ${\bf R}^n$ with compact
support, and then by left or right Clifford-holomorphic functions on
neighborhoods of $A$ using the formula.  This does not work when $A$
has positive Lebesgue measure and finite perimeter, because of
conditions on the distributional derivatives of the product of the
characteristic function of $A$ and functions that are Clifford
holomorphic on $A$.  A classical version of this in one complex
variable is described in Exercise 2 at the end of Chapter 20 of
\cite{r}.

        If $A$ is totally disconnected, then every continuous function
on $A$ can be approximated uniformly by locally constant functions.
In particular, there are ``fat Cantor sets'' which are totally
disconnected and have positive Lebesgue measure.

        For any closed set $A \subseteq {\bf R}^n$, a
continuously-differentiable function on $A$ in the sense of the
Whitney extension theorem is a continuous function $f$ on $A$ together
with a continuous differential $df_x$, $x \in A$, with the same
properties as an ordinary $C^1$ function.  In the present setting, $f$
would take values in $\mathcal{C}(n)$, and $df_x$ would be a
real-linear mapping from ${\bf R}^n$ into $\mathcal{C}(n)$.  The
Whitney extension theorem says that there is an ordinary $C^1$
function on ${\bf R}^n$ with values in $\mathcal{C}(n)$ equal to $f$
on $A$ and whose differential is equal to $df_x$ for $x \in A$.

        If $\Lambda$ is a real-linear mapping from ${\bf R}^n$ into
$\mathcal{C}(n)$, then $\mathcal{D}_L \Lambda, \mathcal{D}_R \Lambda
\in \mathcal{C}(n)$ can be defined in the usual way.  If $f$ is a
continuously-differentiable function on $A$ with values in
$\mathcal{C}(n)$, then this can be applied to $df_x$ to get
$\mathcal{D}_L f(x)$ and $\mathcal{D}_R f(x)$ for each $x \in A$.
This is equivalent to applying $\mathcal{D}_L$ or $\mathcal{D}_R$ to a
$C^1$ extension of $f$ to ${\bf R}^n$ whose differential at $x \in A$
is $df_x$.

        In particular, $\mathcal{D}_L f = 0$ and $\mathcal{D}_R f = 0$
make sense on $A$.  In terms of $C^1$ extensions on ${\bf R}^n$, this
means that $\mathcal{D}_L$ or $\mathcal{D}_R$ of the extension tends
to $0$ as a point approaches $A$, respectively.  More precise rates of
vanishing near $A$ can be obtained from stronger regularity.

         It is easy to see that a real-linear mapping $\Lambda : {\bf
R}^n \to \mathcal{C}(n)$ is uniquely determined by its restriction to
any hyperplane in ${\bf R}^n$ when $\mathcal{D}_L \Lambda = 0$ or
$\mathcal{D}_R \Lambda = 0$.  Conversely, if $\lambda$ is a
real-linear mapping from a hyperplane in ${\bf R}^n$ into
$\mathcal{C}(n)$, then there are extensions of $\lambda$ to
real-linear mappings $\Lambda_1, \Lambda_2 : {\bf R}^n \to
\mathcal{C}(n)$ such that $\mathcal{D}_L \Lambda_1 = 0$,
$\mathcal{D}_R \Lambda_2 = 0$.

        If $A$ is contained in a nice $C^1$ hypersurface in ${\bf
R}^n$ and $f : A \to \mathcal{C}(n)$ is continuously-differentiable,
then $df_x$ is not uniquely determined by $f$ on $A$.  For each $x \in
A$, $df_x$ can be changed in the normal direction to the tangent space
of the hypersurface at $x$.  Thus one can arrange to have
$\mathcal{D}_L f = 0$ or $\mathcal{D}_R f = 0$ on $A$ in this case.
By contrast, there are plenty of fractal sets which are not flat in
this way, so that the differential is uniquely determined by the
function on $A$.  Therefore $\mathcal{D}_L f = 0$ or $\mathcal{D}_R f
= 0$ can be a significant condition on a fractal set.

        On some fractal sets $A$, there may be a lot of nonconstant
continuously-differentiable functions $f$ for which $df_x = 0$ for
each $x \in A$, and hence $\mathcal{D}_L f = \mathcal{D}_R f = 0$.
This includes locally constant functions on totally disconnected sets,
as well as some functions on connected snowflake sets.  If every pair
of elements of $A$ can be connected by a rectifiable curve, then any
continuously-differentiable function $f$ on $A$ with $df_x = 0$ for
each $x \in A$ is constant.  Sierpinski gaskets and carpets and Menger
sponges are examples of such sets.  These sets need not be flat in
${\bf R}^n$, so that the differential is determined by the function on
$A$.

        Even in this type of restricted situation, there are still the
usual problems with products of functions and their $\mathcal{D}_L$ or
$\mathcal{D}_R$ derivatives on $A$.  If there are nonconstant
functions $f$ with $df = 0$ on $A$, then they can be treated much like
constants.

        In one complex variable, holomorphicity can be characterized
by expressing the differential of a function by multiplication in the
complex numbers.  This also makes sense in higher dimensions, using
left or right multiplication in the quaternions or a Clifford algebra,
and is well known to be too restrictive in that the solutions are
affine.  For continuously-differentiable functions on a closed set $A
\subseteq {\bf R}^n$, there are more possibilities.  Functions with
vanishing differentials trivially have this property, and as an
expansion of affine functions one can allow the coefficients to be
continuously-differentiable functions with vanishing differentials.
If $A$ is contained in a $C^1$ curve, then one can choose the
differentials to be given by left or right multiplication in
$\mathcal{C}(n)$.

        Suppose that $A$ is a chord-arc curve in ${\bf R}^n$, so that
the length of an arc on $A$ is bounded by a constant multiple of the
distance between its endpoints.  For example, $A$ might be the graph
of a Lipschitz mapping on a line.  One can check that there is a
continuously-differentiable function $f$ on $A$ with any prescribed
continuous family of differentials $df_x$.  The function is affine
when the differentials are constant, and continuous differentials lead
to approximately affine functions on small arcs.  In particular,
$df_x$ may be defined by multiplication on the left or right by a
continuous $\mathcal{C}(n)$-valued function on $A$.

        Let us identify $x = (x_1, \ldots, x_n) \in {\bf R}^n$ with
$x_1 \, e_1 + \cdots + x_n \, e_n \in \mathcal{C}(n)$. In
$\mathcal{C}(n)$,
\begin{equation}
        x \, x = -(x_1^2 + \cdots + x_n^2) = - |x|^2.
\end{equation}
If $x \ne 0$, then it follows that
\begin{equation}
        x \, \Big(\frac{-x}{|x|^2}\Big) = \Big(\frac{-x}{|x|^2}\Big) \, x = 1.
\end{equation}
Thus $x$ is invertible in $\mathcal{C}(n)$.  There are variants of
this for quaternions, and in ${\bf R}^{n + 1}$ where the additional
coordinate corresponds to the copy of ${\bf R}$ in $\mathcal{C}(n)$.

        Let $b$ be a $\mathcal{C}(n)$-valued function on a closed set
$A \subseteq {\bf R}^n$, and consider
\begin{equation}
        (b(x) - b(y)) (x - y)^{-1} \hbox{ or } (x - y)^{-1} (b(x) - b(y))
\end{equation}
for $x, y \in A$, $x \ne y$.  If $b$ is continuously-differentiable on
$A$ with differential given by left or right multiplication in
$\mathcal{C}(n)$, then these expressions can be extended continuously
to $x = y$, respectively.  These expressions are very similar to
kernels of commutator operators, but one should be careful about
noncommutativity.  They certainly are kernels of commutator operators
when $b$ is real-valued, and can be viewed as linear combinations of
commutators otherwise.

        Analogous expressions in the complex plane are very pleasant.
If $b(z)$ is a polynomial or a rational function on ${\bf C}$, then
$(b(z) - b(w)) / (z - w)$ is a finite sum of products of functions of
$z$ and $w$ individually, which corresponds to an operator of finite
rank.  If $b(z)$ is holomorphic on some open set, then the ratio
extends holomorphically across $z = w$.  If $b(x) = \alpha + \beta \,
x$ or $\alpha + x \, \beta$ on ${\bf R}^n$ for some $\alpha, \beta \in
\mathcal{C}(n)$, then the appropriate quotient is equal to $\beta$,
but noncommutativity prevents one from going further.  For smooth
functions $b$ on smooth curves in ${\bf R}^n$, one gets smoothness of
the quotient in terms of the parameterization of the curve, but this
is not quite the same thing.

\end{document}